\documentclass[12pt]{article}
\usepackage[dvips]{graphicx}
\usepackage{latexsym}
\usepackage{amssymb}

\newcommand{\C}{\mathbb{C}}
\newcommand{\CP}{\mathbb{CP}}
\newcommand{\RP}{\mathbb{RP}}
\newcommand{\R}{\mathbb{R}}
\newcommand{\Z}{\mathbb{Z}}

\newcommand{\PP}{\mathbb{P}}

\renewcommand{\d}{\mathrm{d}}

\newcommand{\koniec}{\begin{flushright}  $\Box $ \end{flushright}}
\def\be{\begin{equation}}
\def\ee{\end{equation}}

\def\Sm{\Sigma}

\def\Om{\Omega}
\def\Th{\Theta}

\def\O{\cal O}
\def\om{\omega}

\def\ov{\overline}

\def\p{\partial}

\def\ov{\overline}

\def\K{\kappa}

\def\a{\alpha}
\def\l{\lambda}

\def\O{{\cal O}}

\newtheorem{theo}{Theorem}[section] 
\newtheorem{prop}[theo]{Proposition}  
\newtheorem{lemma}[theo]{Lemma}
\newtheorem{defi}[theo]{Definition}
\newtheorem{col}[theo]{Corollary}

\begin{document}
\title{Anti-self-dual four-manifolds 
with a parallel real spinor.}

\author{Maciej  Dunajski\thanks{email: dunajski@maths.ox.ac.uk}
\\ The Mathematical Institute,
24-29 St Giles, Oxford OX1 3LB, UK}  

\date{} 
\maketitle
\abstract{ Anti-self-dual metrics in the $(++--)$ signature which admit a covariantly constant real spinor are studied. It is shown that finding such metrics reduces to solving a fourth order integrable PDE, and some examples are given.
The corresponding twistor space is characterised by existence of 
a preferred non-zero real section of $\kappa^{-1/4}$, where $\kappa$
is the canonical line bundle of the twistor space.
It is demonstrated that
if the parallel spinor is preserved by a Killing vector, then the fourth 
order PDE reduces to the
dispersionless Kadomtsev--Petviashvili equation
and its linearisation. Einstein--Weyl structures
on the space of trajectories of the symmetry are characterised by the
existence of a parallel weighted null vector.
} 
\noindent
\section{Introduction}
Constraints on a (pseudo) Riemannian geometry imposed 
by the existence of a parallel spinor essentially depend on 
the properties of the Clifford algebra and the spin group associated with
the metric.  
There has been an interest in such geometries in pure mathematics,
because they extend a list of holonomy groups of 
Riemannian manifolds (where the existence of a parallel spinor implies
Ricci flatness). 
Pseudo Riemannian metrics in various dimensions 
with a covariantly constant spinor have also attracted a lot of attention 
in physics, as such spinors play 
a central role in supersymmetry.

In \cite{B00} Bryant analysed 
all cases  up to six dimensions, together with some higher-dimensional  
examples  of  particular interest. In this paper I shall concentrate on the
four-dimensional case.

Let $({\cal M}, g)$ be a (pseudo) Riemannian spin four-manifold.
Therefore there exist complex two-dimensional vector bundles 
$S_\pm$  (spin-bundles) equipped with parallel symplectic structures 
$\varepsilon_\pm$ such that
\begin{itemize}
\item $\C\otimes T{\cal M}\cong {S_+}\otimes {S_-}$
is a  canonical bundle isomorphism.
\item $g(v_1\otimes w_1,v_2\otimes w_2)
=\varepsilon_+(v_1,v_2)\varepsilon_-(w_1,w_2) $
for $v_1, v_2\in \Gamma(S_+)$ and $w_1, w_2\in \Gamma(S_-)$.
\end{itemize}
I shall  assume that 
there exists a spinor $\iota$ parallel with respect to a Levi-Civita
connection $\nabla$ of $g$
\[
\iota=(p, q)\in\Gamma(S_+), \qquad \nabla \iota=0.
\]
There are three
possible situations depending on the signature of the metric.

\begin{itemize}
\item In the Lorentzian signature $(+---)$
\[
\mbox{Spin(3, 1)}=SL(2,\C),\qquad \iota=(p, q)\rightarrow
\overline{\iota}=(\overline{p}, \overline{q})\in\Gamma(S_-),\qquad
\nabla\overline{\iota}=0.
\]
Therefore $l=\iota\otimes\overline{\iota}$ is a parallel null
vector. 
This condition has been extensively studied in general 
relativity \cite{KSM80}: There exist two real functions $u, v$, 
and one complex function $\xi$ such that
\[
g=\d u\d v-\d \xi\d\overline{\xi}+H(u, \xi,\overline{\xi})\d u^2.
\]
The Ricci flat condition implies that 
$H(u, \xi,\overline{\xi})=\mbox{Re}(f(u, \xi))$ where $f$ is
holomorphic in $\xi$. These solutions are known as {\em p.p waves}.
Analysis of curvature invariants shows that the function
$H(u, \xi,\overline{\xi})$ cannot be eliminated by a 
coordinate transformation.
\item In the Euclidean signature $(++++)$ 
\[
\mbox{Spin(4, 0)}=SU(2)\times\widetilde{SU}(2),
\qquad \iota=(p, q)\rightarrow
\overline{\iota}=(\overline{q}, -\overline{p})\in\Gamma(S_+),\qquad
\nabla\overline{\iota}=0.
\]
A spinor and its complex conjugate form a basis of a spin space $S_+$.
A four-dimensional Riemannian manifold which admits a covariantly
constant spinor must therefore  be hyper-K{\"a}hler.
Hyper K{\"a}hler four-manifolds have been much studied for the last twenty 
five years.
See \cite{Da00} and references therein.
\item In the split signature $(++--)$
(also called
ultra-hyperbolic, Kleinian or neutral)
\[
\mbox{Spin(2, 2)}=SL(2,\R)\times\widetilde{SL}(2, \R),
\qquad \iota=(p, q)\rightarrow
\overline{\iota}=(\overline{p}, \overline{q})\in\Gamma(S_+),\qquad
\nabla\overline{\iota}=0,
\]
and the representation space of the spin group  splits into
a direct sum of two real two-dimensional spin spaces $S_+$ and $S_-$.
The conjugation of spinors is involutive and maps each spin space  
onto itself, and there exists an invariant notion of {\it real spinors}.
\end{itemize}
One can therefore  look for $(++--)$ metrics with a
parallel real spinor (which we choose to be
$\iota\in\Gamma(S_+)$). These metrics do not have to be Ricci-flat.
The resulting geometry will be studied in the rest of this paper.
The isomorphism ${\Lambda^2}_+({\cal M})\cong {\mbox{Sym}}^2(S_+)$ between 
the bundle of self-dual two-forms and the symmetric tensor product of 
two spin bundles implies that the real self-dual two-form 
$\Sm=\iota\otimes\iota\otimes\varepsilon_+$
is covariantly constant
and null (i.e. $\Sm\wedge\Sm=0$) which motivates the following definition:
\begin{defi}
A null-K\"ahler structure on a four-manifold consists of metric 
of signature $(++--)$ and a real spinor field parallel with 
respect to this inner product.
A null-K{\"a}hler structure is anti-self-dual (ASD)
if the self-dual part of the Weyl spinor vanishes.
\end{defi}
The ASD condition on null-K\"ahler structures is worth studying for
at least two reasons:
firstly, four-dimensional vacuum 
metrics in signature $(++--)$ appeared in describing the bosonic sector
of the $N=2$ super-string \cite{OV90}, \cite{BG94}.
Secondly, many integrable systems in dimensions two and three
arise (together with their twistor description) as symmetry
reductions of anti-self-duality equations on $(++--)$ 
background \cite{Wa85,MW96}.

The first aspect will be not 
discussed in the present paper, but I shall reveal some connections
with integrable systems in Sections \ref{Sisystem} and \ref{Ssymmetry}.

In the next Section the ASD null K\"ahler condition
will be related to Einstein-Maxwell equations.
 In Section \ref{Sisystem}
ASD null-K\"ahler condition will be reduced 
to a single fourth order integrable PDE.
Explicit solutions to this equations will provide some examples of
ASD  null-K\"ahler structures.
The resulting twistor theory 
will be described in Section \ref{Sttheory}. The existence of a parallel real spinor 
will be characterised by a real structure preserving a
preferred non-zero section of $\kappa^{-1/4}$, where $\kappa$
is the canonical line bundle of the twistor space.
In Section \ref{Ssymmetry} it will be shown that ASD  null K\"ahler structures
with a symmetry which preserves the parallel spinor are locally given by
solutions to the dispersionless Kadomtsev--Petviashvili equation
and its linearisations \cite{D00}. Einstein--Weyl structures
on the space of trajectories of the symmetry will be characterised by the
existence of a parallel weighted null vector.
The two-component spinor notation will used in the paper. 
The spin spaces $S_-$, and $S_+$ will be denoted by 
$S^{A}$ and $S^{A'}$ respectively.
From now on
the parallel real spinor 
and the corresponding null-K\"ahler two-form 
will be denoted by  
$\iota^{A'}\in\Gamma(S^{A'})$  and 
$\Sm^{0'0'}=\iota_{A'}\iota_{B'}\Sm^{A'B'}\in {\Lambda^2}_+({\cal
M})$. 
The notation is summarised in the Appendix. 
\section{Null-K\"ahler metrics in four dimensions} 
It is well know \cite{Po92} that K\"ahler four-manifolds
with vanishing scalar curvature are necessary ASD. 
This is not true for
scalar-flat 
 null-K\"ahler four-manifolds. Instead one has the following result:
\begin{prop}
\label{ASD-scalar}
Let $\iota_{A'}$ be a parallel real spinor on 
dimensional,  ultrahyperbolic manifold.
Then  the scalar curvature vanishes, the Ricci tensor is null, 
and the self-dual Weyl
spinor is given by
\be
\label{nullC}
C_{A'B'C'D'}=c\iota_{A'}\iota_{B'}\iota_{C'}\iota_{D'}
\ee
for some function $c$ such that $\iota^{A'}\nabla_{AA'}c=0$.
\end{prop}
{\bf Proof:} Vanishing of the scalar curvature follows from the second 
Ricci identity (\ref{ricci2}).
The first  identity (\ref{ricci1}) implies that 
$\Phi_{ABA'B'}\iota^{B'}=0$, so
$\Phi_{ABA'B'}=F_{AB}\iota_{A'}\iota_{B'}$
for some $F_{AB}$.
The formula \ref{nullC} is a direct consequence of (\ref{ricci2}) and  
(\ref{bianchi}) applied to a
covariantly constant spinor.
\koniec
Now I shall show that ASD null-K{\"a}hler metrics can be viewed as
solutions to Einstein equations with electromagnetic stress-energy
tensor. The ASD part of the Maxwell field is given by the Ricci
form, and the SD part is given by the null-K\"ahler form.
\begin{prop}
There is a one-to-one correspondence between
ASD metrics with a constant real spinor, and ASD Einstein-Maxwell
spaces for which the SD part of Maxwell field is null and covariantly
constant.
\end{prop}
{\bf Proof:}
Let $({\cal M}, g)$ be an ASD manifold with a covariantly constant
spinor $\iota_{A'}$. Proposition \ref{ASD-scalar} implies that
\be
\label{Einstein-Maxwell}
\Phi_{ABA'B'}=F_{AB}\iota_{A'}\iota_{B'},
\ee
and the spinor Bianchi identities (\ref{bianchi}) yield
$\nabla^{AA'}F_{AB}=0$. Therefore 
\be
\label{Maxwell_field}
F_{ab}=F_{AB}\varepsilon_{A'B'}+\iota_{A'}\iota_{B'}\varepsilon_{AB}.
\ee
is a Maxwell field, and the formula  (\ref{Einstein-Maxwell}) can be read
off as Einstein equations with a Maxwell 
stress energy tensor 
\[
T_{ab}=\frac{1}{2}\Big(\frac{1}{4}g_{ab}F_{cd}F^{cd}-F_{ac}{F_b}^c\Big)
=F_{AB}\iota_{A'}\iota_{B'}.
\]
Conversely, consider a Maxwell field $F_{ab}$ (\ref{Maxwell_field})
on an ASD background, such that its SD  part is null and constant.
The Maxwell equations give
\[
\nabla^{A}_{A'}F_{AB}=0,\qquad \nabla_{AA'}\iota_{B'}=0,
\]
and the Einstein equations  with the  
stress energy tensor $T_{ab}=F_{AB}\iota_{A'}\iota_{B'}$ yield
(\ref{Einstein-Maxwell}).
\koniec
The following result will be used in Sections \ref{Sisystem}
and \ref{Sttheory}
\begin{prop}
\label{constant_forms}
Let $\Sm^{A'B'}=(\Sm^{0'0'}, \Sm^{0'1'}, \Sm^{1'1'})$ 
be a basis of normalised real SD 
two-forms on an ASD scalar-flat manifold such that
\be
\label{cforms}
\d (\iota_{A'}\iota_{B'}\Sm^{A'B'})=\d (o_{A'}\iota_{B'}\Sm^{A'B'})=0.
\ee
Then there exists a covariantly constant real section of $S_{A'}$.
Conversely, let $\iota_{A'}$ be a covariantly constant section of
$S_{A'}$ on an ASD four manifold. 
Then it is possible to find another section $o_{A'}$, such that
 $(o_{A'}, \iota_{A'})$ form a normalised spin frame, and equations
{\em(\ref{cforms})\em} hold.
\end{prop}
{\bf Proof:} Let $(o_{A'}, \iota_{A'})$ be a normalised spin basis.
The covariant derivatives of the basis can be expressed as 
\[
\nabla_{a}\iota_{B'}=U_a\iota_{B'}+V_ao_{B'},\qquad 
\nabla_{a}o_{B'}=W_a\iota_{B'}-U_ao_{B'}.
\]
The first condition in (\ref{cforms})
can be rewritten as 
${\nabla_A}^{A'}(\iota_{A'}\iota_{B'})=0$, which implies 
\[V_a=2U_{AB'}\iota^{B'}\iota_{A'}.\] The second condition in 
(\ref{cforms}) yields 
$U_{AA'}=\alpha_{A}\iota_{A'}, W_{AA'}=\beta_{A}\iota_{A'}$ for some
$\alpha_{A}, \beta_{A}$. Therefore
\be
\label{cforms1}
\nabla_{AA'}\iota_{B'}=\alpha_{A}\iota_{A'}\iota_{B'}.
\ee
Contracting the RHS of the above equation with ${\nabla^A}_{C'}$,
and symmetrising over $(A'B')$ gives $0$ because $g$ is $ASD$ and
scalar-flat. As a consequence 
${\nabla^{A}}_{(C'}[\iota_{A')}\iota_{B'}\alpha_{A}]=0$ which
gives ${\nabla^A}_{A'}\alpha_{A}=0$. 
Consider the real spinor  
$\hat{\iota}_{A'}:=\iota_{A'}\mbox{exp}f$, where
$\iota^{A'}\nabla_{AA'}f=0$ in order to preserve $\d \Sm^{0'0'}=0$.
Integrability conditions for $\nabla_{AA'}f=\alpha_{A}\iota_{A'}$
are satisfied, as $\alpha_{A}$ solves the neutrino equation.
Therefore we can find  $f$ for each $\alpha_{A}$, and
equation (\ref{cforms1}) implies that $\hat{\iota}_{A'}$ is
covariantly constant.\\
{\bf Converse.}
Now assume $\nabla_{AA'}\iota_{B'}=0$. Consider $\hat{o}_{A'}$ 
such that $\hat{o}_{A'}\iota^{A'}=1$.
The normalization condition
implies $\nabla_{AA'}\hat{o}_{B'}=\gamma_{A}\iota_{A'}\iota_{B'}
+\rho_{A}\hat{o}_{A'}\iota_{B'}$ for some $\gamma_{A}, \rho_{A}$. 
Contracting with ${\nabla^A}_{C'}$ and using Ricci identities (\ref{ricci1})
yields ${\nabla^{A}}_{A'}\rho_{A}=0$, and so 
$\rho_{A}=\iota^{A'}\nabla_{AA'}\phi$. Therefore in the 
null rotated normalised 
spin frame $\iota_{A'}, \hat{o}_{A'}-\phi\iota_{A'}$ we have
$\rho_{A}=0$.
Therefore
\[
{\nabla_{A}}^{A'}(o_{A'}\iota_{B'}+o_{B'}\iota_{A'})=
2{W_A}^{A'}\iota_{A'}\iota_{B'}=\rho_{A}\iota_{B'}=0,
\]
and $\d \Sm^{0'1'}=0$.
\koniec
\section{Anti-Self-Dual null-K\"ahler condition as an integrable system}
\label{Sisystem}
I shall now construct a local coordinate system adapted to the 
parallel spinor, and reduce the ASD null-K\"ahler condition to
a pair of coupled PDEs.
Integrability of these PDEs will be established by using the Lax formulation
(i.e. showing that they arise as the integrability conditions to an over-determined system of linear equations).

Let $S^{A'}={\cal M}\times \C^2$ be the bundle of complex primed spinors.
The natural context for introducing the Lax pair is the geometry of 
the projective primed spin bundle
(also called the correspondence space) 
${\cal F}=\PP(S^{A'})={\cal M}\times\CP^1$.
It is
coordinatized by $(x^{a},\lambda)$, where $x^a$ denotes the coordinates on
$\cal M$ and $\l$ is the coordinate on $\CP^1$ that parametrises the
$\a$-surfaces through $x$ in $\cal M$.  
We relate the fibre coordinates $\pi^{A'}$ on $S^{A'}$ 
to $\lambda$ by
$\lambda=\pi_{0'}/\pi_{1'}$.

Let  $\nabla_{AA'}$ be  a null tetrad of vector fields for the metric $g$ 
on ${\cal M}$ and let $\Gamma_{AA'B'C'}$ be the components of the spin
connection in the associated spin frame.
A  horizontal lift of $\nabla_{AA'}$ to $S_{A'}$ 
given by 
\[
\widetilde{\nabla}_{AA'}=\nabla_{AA'}+\Gamma_{AA'B'C'}\pi^{B'}\frac{\p}{\p\pi_{C'}}
\]
Its horizontality implies $\widetilde{\nabla}_{AA'}\pi_{B'}=0$. 

The space ${\cal F}$ possesses a natural two-dimensional distribution
called the twistor distribution, or Lax pair to emphasise the
analogy with integrable systems.
The Lax pair arises as the
image under the projection $TS^{A'}\longrightarrow T{\cal F}$ of the
distribution spanned by $L_A=\pi^{A'}\widetilde{\nabla}_{AA'}$,
and is given by

\be
\label{laxpair}
L_{0}=\nabla_{00'}-\l \nabla_{01'}+l_{0}\p_{\l},\qquad
L_{0}=\nabla_{10'}-\l \nabla_{11'}+l_{1}\p_{\l}, 
\ee
where
$l_A=\Gamma_{AA'B'C'}\pi^{A'}\pi^{B'}\pi^{C'}$ are cubic polynomials
in $\l$ (note that $\pi_{1'}=1$ in these formulae). 
\begin{theo}{\em\cite{Pe76}}
\label{penroseth}
The twistor distribution on $S_{A'}$ given by {\em(\ref{laxpair})}
is integrable if and only if the Weyl curvature of $g$ is ASD, 
i.e. $C_{A'B'C'D'}=0$.
\end{theo}
We are now ready to reformulate the ASD null-K\"ahler condition
as an integrable system.
\begin{theo}
Real coordinates $(w, z, x, y)$ can be chosen such that all 
ASD null-K\"ahler metric are locally given by
\be
\label{nkmetric}
g=\d w\d x+\d z \d y-\Th_{xx}\d z^2-\Th_{yy}\d w^2+2\Th_{xy}\d w\d z
\ee
where $\Th(w, z, x, y)$ is a solution to a 4th order PDE
$($which we write as a  system of two second order PDEs $)$:
\be
\label{nk1}
\Th_{wx}+\Th_{zy}+\Th_{xx}\Th_{yy}-\Th_{xy}^2=f,
\ee
\be
\label{nk2}
\square f=f_{xw}+f_{yz}+
\Th_{yy}f_{xx}+\Th_{xx}f_{yy}-2\Th_{xy}f_{xy}=0.
\ee
Moreover {\em(\ref{nk1},\ref{nk2})\em} arise as an integrability condition
for the linear system $L_0\Psi=L_1\Psi=0$, where
$\Psi=\Psi(w, z, x, y, \l)$
and
\begin{eqnarray}
\label{Lax2}
L_0&=&(\p_w-\Th_{xy}\p_y+\Th_{yy} \p_x)-\l\p_y+f_y\p_{\l},\nonumber\\
L_1&=&(\p_z+\Th_{xx}\p_y-\Th_{xy} \p_x)+\l\p_x-f_x\p_{\l}.
\end{eqnarray}
\end{theo}
{\bf Proof:}
Let $e^{AA'}$ be a tetrad of real independent one-forms. 
The parallel spinor $\iota_{A'}$ enables us to choose coordinates
$w^{A}=(w, z)$ such that $e^{A0'}=-\iota_{A'}e^{AA'}=\d w^{A}$.
Proposition \ref{constant_forms} implies that we can choose $o_{A'}$
such that $o_{A'}\iota^{A'}=1$, and
\[
\Sm^{0'1'}=\frac{1}{2}\varepsilon_{AB}o_{A'}\iota_{B'}e^{AA'}\wedge e^{BB'}
=o_{A'}e^{AA'}\wedge\d w_{A}
\]
is a closed two-form. Therefore the Frobenius theorem
guarantees the existence of coordinates $x_{A}=(x, y)$ such that
${e_A}^{1'}=o_{A'}{e_A}^{A'}=\d x_{A}+\Th_{AB}\d w^{B}$, where
$\Th_{AB}=\Th_{AB}(w, z, x, y)$ is symmetric in $A$ and $B$.
With this choice $\Sm^{0'1'}=\d x_A\wedge \d w^A$, and the metric is given by
\[
g=\d x_A\d w^A+\Th_{AB}\d w^A\d w^B.
\]
Calculating the components of the spin connection yields
$\Gamma_{AA'B'C'}=A_{AA'}\iota_{B'}\iota_{C'}$. The
residual conformal freedom is used to set  $A_{AA'}=\beta_{A}\iota_{A'}$.

The tetrad of vector fields $\nabla_{AA'}$ dual to $e^{AA'}$ is 
\[
\nabla_{A1'}=\iota^{A'}\nabla_{AA'}=\frac{\p}{\p
x^{A}},\qquad
\nabla_{A0'}=o^{A'}\nabla_{AA'}=\frac{\p}{\p w^A}+\Th_{AB}\frac{\p}{\p x_B},
\]
and the Lax pair (\ref{laxpair}) is  
\[
L_A=\frac{\p}{\p
x^{A}}-\l\Big(\frac{\p}{\p w^A}+\Th_{AB}\frac{\p}{\p x_B} \Big)
+l_A\frac{\p}{\p \l},
\]
where in the chosen spin frame  $l_{A}=\beta_{A}$ do not depend on $\l$.
Consider the Lie bracket
\[
[L_0, L_1]=\Big(\frac{\p\Th^{CD}}{\p w^C}+
\Th_{AC}\frac{\p\Th^{AD}}{\p x_C}
-\beta^D- \l\frac{\p \Th^{AD}}{\p x^A}\Big)\frac{\p}{\p x^D}
+\Big(\frac{\p \beta^A}{\p w^A}+\Th_{AB}\frac{\p \beta^A}{\p x_B} 
+\l\frac{\p \beta_{A}}{\p x_A}\Big)\frac{\p}{\p \lambda}.
\]
The ASD condition is equivalent to integrability of the distribution
$L_A$. In fact $[L_0, L_1]=0$ because there is no $\p/\p\om^A$ term 
in the Lie bracket above. 
We deduce that 
\[
\Th_{AB}=\delta_{A}\delta_{B}\Th, \qquad
\beta_A=\delta_A f,\qquad{\mbox{where}}\;\delta_A:=\iota^{A'}\nabla_{AA'}=
\frac{\p}{\p x^A},
\]
and
$f=f(w, z, x, y)$ and $\Th(w, z, x, y)$ satisfy
\[
\Th_{wx}+\Th_{zy}+\Th_{xx}\Th_{yy}-\Th_{xy}^2=f+{\cal F}(w, z), \qquad
\square f=0.
\]
To obtain (\ref{nk1}) we absorb ${\cal F}(w, z)$ into $f$ without
changing (\ref{nk2}).
\koniec
{\bf Remarks}
\begin{itemize}
\item One-forms $(e^{00'}, e^{10'})$ span a differential ideal,
and $\Sm^{0'0'}=e^{00'}\wedge e^{10'}=\d w\wedge\d z$
is the null-K\"ahler form (it is covariantly constant with respect to
the metric (\ref{nkmetric})).
On the other hand $(e^{01'}, e^{11'})$ do not span an ideal
unless $f=0$ in which case $g$ is pseudo-hyper-K{\"a}hler.
\item
The components of the Ricci and Weyl curvatures, and the Levi-Civita spinor 
connection are 
\begin{eqnarray*}
C_{ABCD}&=&\delta_A\delta_B\delta_C\delta_D\Th,\qquad
C_{A'B'C'D'}=\iota_{A'}\iota_{B'}\iota_{C'}\iota_{D'}\square f=0,\qquad
R=0,\\
\Phi_{ABA'B'}&=&\iota_{A'}\iota_{B'}\delta_{A}\delta_{B}f,\qquad
\Gamma_{AB}=\delta_{A}\delta_{B}\delta_{C}\Th\d w^{C},
\qquad \Gamma_{A'B'}=\iota_{A'}\iota_{B'}\delta_Af \d w^{A}.
\end{eqnarray*}
The Bianchi identity (\ref{bianchi}) is satisfied as a consequence of 
equation (\ref{nk2}).
To sum out $f$ is a potential for a null Maxwell field (so
called Hertz potential),
and $\Th$ is a non-linear potential for a metric.

Potential forms of complexified null Einstein-Maxwell 
equations were given in
\cite{GPR77} and \cite{R00}. It will be instructive to look for a
non-trivial overlap between them and the one given above.
\end{itemize}
\subsection{Examples}
{\bf (1)} Consider a class of metrics given by $\Th_x=0$. Equations
(\ref{nk1},\ref{nk2}) reduce to 
\[
\Th_{yz}=f,\qquad f_{yz}=0.
\]
The general solution
is given by
\be
\Th=B(w,y)+z\int A(w,y)\d y,\qquad f=A(w, y),
\ee
where $A(w, y)$ and $B(w, y)$ are arbitrary functions.
(In fact there are other terms linear in $y$ 
and depending on arbitrary functions of $w, z$. These terms can be
gauged away, as they do not change the metric.)
We have
\[
g=\d w\d x+\d z \d y-(zA_y+B_{yy})\d w^2.
\]
If $A_y=0$ then $g$ is pseudo hyper-K\"ahler.
\\
{\bf (2)} Consider solutions with $f=\Th_{v}$, where $v$ is one of $(w,
z, x, y)$. Equation (\ref{nk2}) implies
\[
\frac{\p}{\p v}(\Th_{wx}+\Th_{zy}+\Th_{xx}\Th_{yy}-\Th_{xy}^2)=0,
\]
and differentiating equation (\ref{nk1}) yields $\Th_{vv}=0$.
It is enough to consider two sub-cases : $\Th_{xx}=0$ and  $\Th_{ww}=0$ 
where the former one can be integrated explicitly
\be
g=\d w\d x+\d z \d y-(xP_y+z(P-P_w+2PP_y)+Q)\d w^2+2P\d w\d z,
\ee
where $P(w, y)$ and $Q(w, y)$ are arbitrary functions.
\\
{\bf (3)} \begin{eqnarray}
\Th&=& A\Big(\frac{x}{y}\Big),\nonumber\\
g&=&\d w\d x+\d z \d y-\frac{1}{y^2}\Big((\frac{x^2}{y^2}A^{''}
+2\frac{x}{y}A^{'})\d w^2+A^{''}\d z^2+2(\frac{x}{y}A^{''}+A^{'})\d w\d z\Big)
\end{eqnarray}
where $A$ is an arbitrary function, and $A^{'}$ is its derivative.\\
{\bf (4)} \begin{eqnarray}
\Th&=& xA(y)+B(y),\nonumber\\
g&=&\d w\d x+\d z \d y-(xA_{yy}+B_{yy})\d w^2+2A_y\d w\d z
\end{eqnarray}
where $A$ and $B$ are arbitrary functions of one variable.
\section{Twistor theory of ASD Null-K\"ahler metrics}
\label{Sttheory}
All ASD null-K{\"a}hler metrics locally arise from solutions to
(\ref{nk1},\ref{nk2}). Non-analytic solutions are generic in $(2,2)$ 
signature. However, in order to find a twistor description, in this section 
I shall restrict myself to real-analytics solutions.

Given an analytic
solution to (\ref{nk1},\ref{nk2}) one can obtain the corresponding
twistor space by equipping ${\cal M}\times \CP^1$ with an integrable complex
structure: The basis of $[0, 1]$ vectors is $(L_0, L_1,
\p_{\ov{\l}})$, where $(L_0, L_1)$ are given by (\ref{Lax2}).
The parallel spinor
$\iota^{A'}$ gives
rise to the section $l=\iota^{A'}\pi_{A'}$ of $\K^{-1/4}$.
In this section I shall perform this 
construction (together with its converse)  in a coordinate independent way.
\begin{defi}
An $\alpha$-surface is a totally null two-dimensional surface, such that a 
two-form orthogonal to its tangent plane is SD.
\end{defi}
There are Frobenius integrability conditions for the existence
of such $\a$-surfaces through each point
and these are equivalent, by Theorem \ref{penroseth} to the vanishing of
the self-dual part of the Weyl curvature, $C_{A'B'C'D'}$.  Thus, given
$C_{A'B'C'D'} =0$, we can define a twistor space ${\cal PT}$ to be the
three complex dimensional manifold of $\a$-surfaces in ${\cal M}$.

A tangent space to an $\alpha$-surface is spanned by null vectors
of the form $\lambda^{A}\pi^{A'}$ with $\pi^{A'}$ fixed and $\lambda^{A}$
arbitrary. As mentioned in the introduction, in the split
signature any spinor has an invariant decomposition into its real and 
imaginary part. A  {\em real $\alpha$-surface} corresponds to both
$\lambda^{A}$ and $\pi^{A'}$ being real.

In general $\pi^{A'}=\mbox{Re}{\pi^{A'}}+i \mbox{Im}{\pi^{A'}}$, and the correspondence 
space ${\cal F}$ defined in the last
section decomposes into two open sets 
\begin{eqnarray*}
{\cal F}_+&=&\{ (x^a, [\pi^{A'}])\in {\cal F};
\mbox{Re}({\pi_{A'}})\mbox{Im}({\pi^{A'}})>0\}={\cal M}\times D_+,\\
{\cal F}_-&=&\{ (x^a, [\pi^{A'}])\in {\cal F};
\mbox{Re}({\pi_{A'}})\mbox{Im}({\pi^{A'}})<0\}={\cal M}\times  D_-,
\end{eqnarray*}
where $D_{\pm}$ are two copies of a Poincare disc.
These complex submanifolds are separated by a real correspondence space
\[
{\cal F}_0=\{ (x^a, [\pi^{A'}])\in {\cal F}; \mbox{Re}({\pi_{A'}})\mbox{Im}({\pi^{A'}})=0\}
={\cal M}\times\RP^1.
\]
The vector fields (\ref{laxpair}) together with the complex structure
on the  $\CP^1$ give ${\cal F}$ a structure of a complex manifold
${\cal PT}$:
The integrable sub-bundle of $T\cal F$ is spanned by $L_0, L_1,  \p_{\ov\l}$.
The distribution  (\ref{laxpair}) with $\l\in\RP^1$
define a foliation of ${\cal F}_0$
with a quotient ${\cal PT}_0$ which leads to 
a double fibration:
\be
\label{doublefib}
{\cal M}\stackrel{p}\longleftarrow 
{\cal F}_0\stackrel{q}\longrightarrow {\cal PT}_0.
\ee
The twistor space  ${\cal PT}$ is a union of two open subsets
${\cal PT}_+=({\cal F}_+)$ and ${\cal PT}_-=({\cal F}_-)$
separated by a three-dimensional real boundary\footnote{ In \cite{Wo92} Woodhouse
performed a careful analysis of the twistor correspondence for flat $(++--)$
metrics, and showed how functions on ${\cal PT}_0=\RP^3$ can be used to
construct smooth solutions to the ultra-hyperbolic wave equation.}
({\em real twistor space}) ${\cal PT}_0:=q({\cal F}_0)$.

The {\em real structure} $\sigma(x^a)=\overline{x}^a$ maps $\alpha$-surfaces
to  $\alpha$-surfaces, and therefore induces an anti-holomorphic
involution $\sigma:  {\cal PT}\rightarrow {\cal PT}$.
The fixed points of this involution  correspond to real  
$\alpha$-surfaces in ${\cal M}$. There is an $\RP^1$ worth 
of such $\alpha$-surfaces through each point of ${\cal M}$.
The set of fixed points of $\sigma$ in  ${\cal PT}$ is 
${\cal PT}_0$.

Each point $x\in\cal M$ determines a sphere $l_x$ made up of all the
$\alpha$-surfaces through $x$. The normal bundle of $l_x$
in $\cal PT$ is $N=T{\cal PT}|_{l_x}/Tl_x$. This is a rank-two vector
bundle over $\CP^1$, therefore it has to be one of the standard 
line\footnote{ Here ${\cal
O}(n)$ denotes the line bundle over $\CP^1$ with transition functions
$\l^{-n}$ from the set $\l\neq\infty$ to $\l\neq 0$ (i.e.\ Chern class
$n$).}
bundles ${\cal O}(n)\oplus {\cal O}(m)$.
\begin{lemma} Let $p:{\cal F}={\cal M}\times\CP^1\longrightarrow {\cal M}$.
The holomorphic curves $l_x:=p^{-1}(x)$, $x\in {\cal M}$ 
have normal bundle $N={\cal
O}(1)\oplus {\cal O}(1)$.   
\end{lemma}
{\bf Proof.} The bundle $N$ can be identified with the
quotient $p^*(T_x{\cal M})/\{ \mathrm{span }\;L_{0},L_{1}\}$.
In their homogeneous form the operators $L_{A}$ have weight one, so the
distribution spanned by them is isomorphic to the bundle
$\C^{2}\otimes{\cal O}(-1)$.  The definition of the normal bundle as
a quotient gives a sequence of sheaves over $\CP^1$.
\[
0\longrightarrow \C^{2}\otimes{\cal O}(-1) \longrightarrow \C^{4}  
\longrightarrow N\longrightarrow 0
\]
and we see that $N={\cal O}(1)\oplus{\cal O}(1) $ , because the
last map, in the spinor notation, 
is given explicitly by $V^{AA'}\mapsto
V^{AA'}\pi_{A'}$  projecting
onto ${\cal O}(1)\oplus{\cal O}(1) $.
\koniec
\smallskip
If ${\cal M}$ is ASD null-K\"ahler then ${\cal PT}$ has an 
additional structure:
\begin{theo}
\label{Twmaci}
Let ${\cal PT}$ be a    
three-dimensional complex manifold with 
\begin{itemize}
\item a four parameter family of
rational curves with normal bundle ${\cal O}(1)\oplus {\cal O}(1)$,
\item a preferred section of $\kappa^{-1/4}$, where $\kappa$ is the 
canonical bundle of ${\cal PT}$,
\item  an
anti-holomorphic involution 
$\rho:{\cal PT}\rightarrow{\cal PT}$ fixing a real equator of  
each rational curve, and leaving the section of
$\kappa^{-1/4}$ above invariant,
\end{itemize}
Then the real moduli space ${\cal M}$ of the $\rho$-invariant 
curves is equipped with 
a restricted conformal class $[g]$ of  ASD null-K\"ahler
metric: If $g\in[g]$ and $\Sm^{0'0'}$ is a null-K\"ahler
two-form
then   
$\hat{g}=\Omega^2 g\in[g]$ for any $\Omega$ such that
$\d\Omega\wedge\Sm^{0'0'}=0$.
Conversely, given a real analytic  ASD null K\"ahler metric, 
there exists a corresponding twistor space with the above
structures.
\end{theo}
{\bf Proof}:
Let $g$ be a real analytic ASD metric with a covariantly constant real
spinor 
$\iota_{A'}$.
From $C_{A'B'C'D'}=0$ it follows that there exist coordinates
$\pi^{A'}$ on the fibres of $S^{A'}\rightarrow{\cal M}$ such that
$
\pi^{A'}\widetilde{\nabla}_{AA'}\pi^{B'}=0.
$
Therefore a parallel section  $\iota_{A'}$ of $S_{A'}$ determines a 
function $l=\pi^{A'}\iota_{A'}$ constant along the twistor distribution.
The line bundle  $\pi^{A'}\iota_{A'}=0$ on ${\cal PT}$
is isomorphic to  $\kappa^{-1/4}$, where $\kappa=\Omega^3{\cal PT}$
is the canonical bundle.\\
{\bf Converse.} The global section $l$ of $\kappa^{-1/4}$, when pulled back to $S_{A'}$ determines a homogeneity degree one function on each fibre of $S_{A'}$ and so must, by globality, be given by $l=\iota^{A'}\pi_{A'}$ and since
$l$ is pulled back from twistor space, it must satisfy 
$\pi^{A'}\widetilde{\nabla}_{AA'}l=0$.
This implies 
\be
\label{twistoreq1}
{\nabla}_{AA'}\iota_{B'}=\varepsilon_{A'B'}\alpha_{A}
\ee
for some $\alpha_{A}$. 
Choose a representative in $[g]$ with $R=0$.
Contracting (\ref{twistoreq1}) with
${\nabla^A}_{C'}$ and using the spinor Bianchi identity gives
\[
{\nabla^A}_{C'}\nabla_{AA'}\iota_{B'}=C_{A'B'C'D'}\iota^{D'}-
\frac{1}{12}R\varepsilon_{C'(B'}\iota_{A')}=0=
\varepsilon_{A'B'}{\nabla^A}_{C'}\alpha_A,
\]
so $\alpha_{A}$ is a solution to the  ASD spin $1/2$  equation 
$\nabla^{AA'}\alpha_A=0$ (so-called neutrino equation).
It can be written in terms of a potential
\be
\label{potentialeq}
\alpha_{A}=\iota^{A'}\nabla_{AA'}\phi
\ee
because the integrability conditions
$\iota^{A'}\iota^{B'}{\nabla^A}_{A'}\alpha_{A}
=\alpha_{A}\iota^{A'}{\nabla^A}_{A'}\iota^{B'}$ are satisfied.
Here $\phi$ is a real analytic function which satisfies
\be
\label{equationforf}
\nabla^a\nabla_a\phi+\nabla_a\phi\nabla^a\phi=0
\ee
as a consequence of the neutrino equation.
Consider a conformal rescaling \[
\hat{g}=\Om^{2}g,\qquad 
\hat{\varepsilon}_{A'B'}=\Omega{\varepsilon}_{A'B'},\qquad 
\hat{\iota}_{A'}=\Omega{\iota}_{A'},\qquad
\hat{\iota}^{A'}={\iota}^{A'},\qquad
{\hat R}=R+\frac{1}{4}\Omega^{-1}\square\Omega.
\]
The twistor equation (\ref{twistoreq1}) is conformally invariant as
${\hat{\nabla}_A}^{(A'}\hat{\iota}^{B')}=
\Om^{-1}{{\nabla}_A}^{(A'}{\iota}^{B')}=0$. 
Choose $\Om\in\ker\square$ so that $\hat{R}=0$. 
Let $\Upsilon_a=\Om^{-1}\nabla_a\Om$. Then
\[
\hat{\nabla}_{AA'}{\hat\iota}^{B'}=
\nabla_{AA'}{\iota}^{B'}+{\varepsilon_{A'}}^{B'}\Upsilon_{AB'}\iota^{C'}
={\varepsilon_{A'}}^{B'}(\iota^{C'}\nabla_{AC'}(\phi+\mbox{ln}\Om))
\]
where we used (\ref{twistoreq1}) and (\ref{potentialeq}).
Notice that, as a consequence of (\ref{equationforf}),
$exp(\phi)\in\ker\square$ and we can choose $\mbox{ln}\Om=-\phi$,
and 
\be
\label{finaleq}
\hat{\nabla}_{AA'}{\hat\iota}^{B'}=0.
\ee
We can still use the residual gauge freedom and add to $\phi$ and
an arbitrary function $\Omega$ constant along $\iota^{A'}\nabla_{AA'}$,
which by Frobenius theorem implies $\d\Omega\wedge\Sm^{0'0'}=0$.
This means (\ref{finaleq}) is invariant under a conformal rescaling
by functions constant along the leaves of the congruence defined by 
$\hat{\iota}^{A'}$. Such conformal transformations do not change
$\hat{R}=0$.
\koniec
\section{ASD null K\"ahler metrics with symmetry}
\label{Ssymmetry}
In this section I shall consider ASD null
K{\"a}hler metrics which admit a Killing vector preserving the
parallel spinor. Let us call them 
{\it ASD null K\"ahler metrics with symmetry}.
I shall show that all such metrics are (at least in the real analytic
case) locally determined by solutions to a certain integrable
equation and its linearisation.

Before establishing  this result I shall review some facts about
Einstein--Weyl (EW) spaces which admit a parallel weighted 
vector\cite{DMT00}.
 \subsection{3D Einstein--Weyl spaces with a parallel weighted vector}
Let ${\cal W}$ be a
three-dimensional real manifold with a torsion-free connection $D$ and a
conformal metric $[h]$.  We shall call ${\cal W}$ a Weyl space if the
null geodesics of $[h]$ are also geodesics for $D$.  This condition is
equivalent to \be
\label{ew1}
Dh=\nu\otimes h
\ee
for some one form $\nu$. Here $h$ is a representative metric in
the conformal class. If we change this representative by
$h\longrightarrow \phi^2 h$, then $\nu\longrightarrow
\nu+2\d\ln{\phi}$.  A tensor object $T$ which transforms as $
T\longrightarrow \phi^m T$ when
$ h\longrightarrow \phi^2 h$
is said to be conformally invariant of weight $m$.
The
covariant derivative of a one-form $\beta$ of weight $m$ 
can be expressed in terms of the Levi-Civita connection of $h$:
\be
\label{wcowder}
\widetilde{D}\beta= \nabla\beta-\frac{1}{2}(
\beta\otimes\nu+(1-m)\nu\otimes\beta-h(\nu, \beta)h).
\ee
The conformally invariant Einstein--Weyl (EW) condition on
$({\cal W}, h, \nu)$ is
\[
W_{(ij)}=\frac{1}{3}\Lambda h_{ij}.
\]
Where $\Lambda$ and $W_{ij}$ are the  scalar curvature, 
and the Ricci tensor of the Weyl connection.

Three-dimensional EW structures are related to four-dimensional
ASD conformal structures by the Jones-Tod correspondence \cite{JT85}:
\begin{prop}[\cite{JT85}]
\label{prop_JT}
Let $({\cal M}, \hat{g})$ be a ultra-hyperbolic 
four-manifold with  ASD conformal curvature,
and a  conformal Killing
vector $K$.
The EW structure in indefinite signature 
on the space  ${\cal W}$ of trajectories of $K$ 
is defined by
\be
\label{EWs}
h:=|K|^{-2}\hat{g}-|K|^{-4}{\bf K}\odot {\bf K},\;\;\; \nu:=2|K|^{-2}
\ast_{\hat{g}}({\bf K}\wedge \d{\bf K}),
\ee
where $|K|^2:=\hat{g}_{ab}K^aK^{b}$, ${\bf K}$ is the one form dual to $K$ and
$\ast_{\hat{g}}$ is taken with respect to $\hat{g}$. 
All three-dimensional EW structures arise in this way.

Conversely, let $(h, \nu)$ be a three--dimensional EW structure 
with a signature $(++-)$ 
on ${\cal W}$, and
let $(V, \a)$ be a pair consisting of a function of weight $-1$ 
and a one-form on ${\cal W}$ which satisfy
the generalised monopole equation
\be
\label{EWmonopole}
\ast_h(\d V+(1/2)\nu V) =\d\a,
\ee
where $\ast_h$ is taken with respect to $h$. Then
\be
\label{Vag}
g=Vh-V^{-1}(\d z+\a)^2
\ee
is an ASD metric with an isometry $K=\p_z$.
\end{prop}
In \cite{DMT00} it has been demonstrated that
if an EW space  admits
a parallel weighted vector, the coordinates can be found in which 
the metric and the one form are given by 
\be
\label{EWdkp}
h=\d y^2-4\d x\d t-4u\d t^2,\qquad\nu=-4u_{x}\d t,\qquad u=u(x, y, t)
\ee
and the EW equations reduce to 
the dispersionless Kadomtsev--Petviashvili equation
\be
\label{dkp}
(u_t-uu_x)_x=u_{yy}.
\ee
If $u(x, y, t)$ is a smooth real function of real variables then
(\ref{EWdkp}) has signature $(++-)$.
It  has also been  shown that 
there exists  a twistor construction of 
EW spaces (\ref{EWdkp}) given by the following Theorem
\begin{theo}{\em\cite{DMT00}}
\label{dKPtwistor}
There is a one to one correspondence 
between Einstein-Weyl spaces ($\ref{EWdkp}$)  
obtained from solutions to the equation $(\ref{dkp})$
and two-dimensional complex manifolds with
\begin{itemize}
\item A three parameter family of rational curves with normal bundle
$\O(2)$.
\item
A global section $l$ of $\K^{-1/4}$, where $\K$ is the canonical bundle.
\item
An anti-holomorphic involution fixing a real slice, leaving a rational
curve and the preferred section of
$\K^{-1/4}$ invariant.
\end{itemize}
\end{theo}
\subsection{Symmetry reduction  }
Now we are ready to establish the main result of this section.
\begin{theo}
\label{null_kahler_symm}
Let $H=H(x, y, t)$ and $W=W(x, y, t)$ be 
smooth real-valued functions on an open set
${\cal W}\subset\R^3$ which satisfy\footnote{
With definition $u=H_x$ the  $x$ derivative of equation 
(\ref{Heqn}) becomes 
the dispersionless Kadomtsev--Petviashvili equation (\ref{dkp}) 
originally used in \cite{DMT00}.
There are some computational advantages in working with the
`potential` form (\ref{Heqn}).}  
\be
\label{Heqn}
H_{yy}-H_{xt} +H_xH_{xx}=0,
\ee
\be
\label{lindKP}
W_{yy}-W_{xt}+(H_xW_x)_x=0.
\ee
Then
\be
\label{scalarflat}
g=W_x(\d y^2-4\d x\d t-4H_x\d t^2)-W_x^{-1}(\d z-W_x\d y-2W_y\d t)^2
\ee
is an  ASD null K\"ahler metric on a circle bundle 
${\cal M}\rightarrow {\cal W}$. All real analytic ASD null  K\"ahler metrics with symmetry
arise from this construction.
\end{theo}
{\bf Proof.}
\label{proof}
Let $(h, \nu)$ be a three--dimensional EW structure given by
(\ref{EWdkp})  (with $u=H_x$)  and let $(V, \a)$ 
be a pair consisting of a function  and a one-form
which satisfy the generalised monopole equation
(\ref{EWmonopole}).
The ultra-hyperbolic metric
\be
\label{Vag1}
g=V(\d y^2-4\d x\d t-4H_x\d t^2)-V^{-1}(\d z+\a)^2
\ee
is therefore  ASD. It satisfies ${\cal L}_Kg=0$, where  $K=\p_z$. 
Using the relations
\[
\ast_h\d t=\d t\wedge\d y,\qquad \ast_h\d y=2\d t\wedge\d x,\qquad 
\ast_h\d x=\d y\wedge\d x+2H_x\d y\wedge\d t
\]
we verify that equation (\ref{lindKP}) is equivalent to 
$\d\ast_h(\d+\nu/2)(W_x)=0$. Therefore 
\[
W_{xx}\d y\wedge\d x+(2(H_xW_x)_x-W_{tx})\d y\wedge\d t+
2W_{xy}\d t\wedge\d x=\d\a,
\]
and we deduce that $V=W_x$ is a  
solution to the monopole equation (\ref{EWmonopole})
on the EW background given by (\ref{EWdkp}).
We choose a gauge in which $\a=Q\d y+P\d t$. This yields
\be
\label{itermidia}
Q_x=-W_{xx},\qquad P_x=-2W_{xy},\qquad P_y-Q_t=2(H_xW_x)_x-W_{xt},
\ee
so $Q=-W_x+A(y,t),\; P=-2W_y+B(y,t)$ and
$
\a=-W_x\d y-2W_y\d t+A\d y+B\d t.
$
The integrability conditions $P_{xy}=P_{yx}$ are given by
(\ref{lindKP}), and $A_t=B_y$. Therefore there exists $C(y, t)$ such
that $A=C_y, B=C_t$.
We now replace $z$  by $z-C$ and the metric
(\ref{Vag1}) becomes (\ref{scalarflat}). This proves that 
(\ref{scalarflat}) is ASD. It is also scalar-flat, 
because, as a consequence of (\ref{lindKP}),
\be
\label{scalar}
R=8(W_{xyy}-W_{xxt}+(H_xW_x)_{xx})W_x=0.
\ee
We now choose the null tetrad 
\begin{eqnarray*}
e^{00'}&=&-2W_x\d t,\;e^{10'}=\frac{\d z-2W_y\d t}{2W_x},\\
e^{01'}&=&\d z-2W_x\d y-2W_y\d x+ze^{00'},\;
e^{11'}=\d x+H_x\d y+ze^{10'},
\end{eqnarray*}
such that 
$g=2(e^{00'}e^{11'}-e^{10'}e^{01'})$. 
The basis of SD two forms $\Sm^{A'B'}$ is given by 
\begin{eqnarray*}
\Sm^{0'0'}&=&=\iota_{A'}\iota_{B'}\Sm^{A'B'}= e^{00'}\wedge e^{10'}=\d z\wedge\d t\\
\Sm^{0'1'}&=&=\iota_{A'}o_{B'}\Sm^{A'B'}= 
e^{10'}\wedge e^{01'}-e^{00'}\wedge e^{11'}=
\d t\wedge \d (z^2)+2\d t\wedge \d W +\d y\wedge\d z.\\
\Sm^{1'1'}&=&=o_{A'}o_{B'}\Sm^{A'B'}= e^{01'}\wedge e^{11'}=2W_x\d x\wedge\d y+
2(zW_x+W_y)\d x\wedge\d t-
\d x\wedge\d z\\
& &+(2H_xW_x-2zW_y)\d t\wedge\d y+z\d z\wedge\d y+(H_x+z^2)\d z\wedge\d t.
\end{eqnarray*}
These two-forms satisfy
\[
-2\Sm^{0'0'}\wedge\Sm^{1'1'}=\Sm^{0'1'}\wedge\Sm^{0'1'},
\qquad \d \Sm^{0'0'}=0,\qquad \d \Sm^{0'1'}=0,
\]
\be
\label{sigma11}
\d \Sm^{1'1'}=
\d(H_{x}-2W)\wedge\d t\wedge\d z
+(W_{xt}-W_{yy}-(H_xW_{x})_x)\d x\wedge\d y\wedge\d t.
\ee
Therefore Proposition \ref{constant_forms} implies that 
the metric
(\ref{scalarflat}) admits a constant spinor
which is preserved by $K=\p_z$.
\\  
{\bf Converse :} 
Let $g$ be a real analytic ASD metric with a covariantly constant spinor 
$\iota_{A'}$, which is Lie derived along a Killing vector $K$.
Theorem \ref{Twmaci} implies that the corresponding twistor space ${\cal
PT}$ is equiped with $l\in\Gamma(\kappa^{-1/4})$
The Killing vector  $K$ gives rise to a holomorphic vector field on ${\cal PT}$
which preserves $l$ . Therefore the minitwistor space ${\cal Z}$
(the space of trajectories of $K$ in $\cal PT$) 
also admits a preffered real section  of  
the $-1/4$ power of its  canonical bundle. 
The minitwistor space $\cal Z$
satisfies the assumptions of Theorem \ref{dKPtwistor} 
and the corresponding EW metric is of the form $\hat{g}=\Omega^2g$,
where $g$  is given by   (\ref{scalarflat}). Both $\hat{g}$ and $g$
are scalar flat (this  follows from the spinor Ricci identities and 
from equation (\ref{scalar}) respectively).
As a consequence  we  deduce that $\Omega=\Omega(t)$.
Now we can use the coordinate freedom \cite{DMT00} to absorb
$\Omega$ in the solution to the  equation (\ref{Heqn}).
\koniec
\begin{col}
Let $({\cal M}, g)$ be an ASD null-K\"ahler manifold with symmetry.
Then the Einstein--Weyl structure induced by {\em(\ref{EWs})} on the space of orbits
of this symmetry is locally of the form {\em(\ref{EWdkp})}.
\end{col}
{\bf Remarks:}
\begin{itemize}
\item 
Theorem \ref{null_kahler_symm} is analogous to a result of LeBrun \cite{L91} who 
constructs all scalar-flat K\"ahler metrics 
with symmetry in Euclidean signature from solutions to the $SU(\infty)$
Toda equation and its linearisation.
\item
A class of solutions to the monopole 
(\ref{lindKP}) can be obtained from vectors tangent to a space of 
solutions to (\ref{Heqn}) generated by
\[
x\rightarrow a(\epsilon)x+b(\epsilon),\qquad
t\rightarrow a(\epsilon)t+c(\epsilon),\qquad 
y\rightarrow a(\epsilon)y+e(\epsilon),
\]
where $a(0)=a, b(0)=b, c(0)=c, e(0)=e$. The corresponding linearised 
solution is given by
\[
W(x, y, t)=\frac{\d H}{\d \epsilon}|_{\epsilon=0}
=a(xH_x+yH_y+tH_t)+bH_x+cH_t+eH_y.
\]
\item
 If $H=const$ then (\ref{lindKP}) reduces to the wave
equation in $2+1$ dimension, and consequently the metric 
(\ref{scalarflat}) is the $(++--)$ Gibbons--Hawking solution \cite{GH78}.
\item Note that $\d \Sm^{1'1'}\neq0$ unless $W=H_x/2+f(t)$, in which case 
\[
\d \Sm^{1'1'}=\d(H_{xt}-H_xH_{xx}-H_{yy})\wedge\d y\wedge\d t=0,
\]
and we are working in a covariantly constant real spin frame.
The metric 
\be
\label{HCdKP1}
g=\frac{H_{xx}}{2}(\d y^2-4\d x\d t-4H_x\d t^2)
-\frac{2}{H_{xx}}(\d z-\frac{H_{xx}\d y}{2}-H_{xy}\d t)^2
\ee
is therefore pseudo 
hyper-K{\"a}hler.  In \cite{DMT00} it was shown that all 
pseudo hyper-K{\"a}hler metrics with a symmetry satisfying
$\d K_+\wedge\d K_+=0$ are locally given by (\ref{HCdKP1}).
Here $\d K_+$ is a self-dual part of $\d K$.
\item
If $W_x\neq H_{xx}/2$ then (\ref{scalarflat}) is not Ricci--flat.
This can be verified by a direct calculation. 
It also follows from
 more geometric reasoning: The  Killing vector 
$K=\p_z$ acts on  SD
two-forms  by a Lie derivative. One can choose a basis $\Sm^{A'B'}$
such that one element of this basis is fixed, and the Killing vector
rotates the other two. The components of the SD derivative of
$K$  are coefficients of these rotations. Therefore 
$(\d K)_+=const$ if $g$ is pseudo hyper-K\"ahler. In our case
$
\d K_+=({H_{xx}}/{W_x})\d z\wedge\d t.
$
Therefore $H_{xx}/W_x$ must be constant for  
(\ref{scalarflat}) to be Ricci-flat.
An example of a non-vacuum metric 
is given by $W=H_{y}/2$.   
\end{itemize}
\subsection{Pseudo hyper-K\"ahler metrics with symmetry}
In this subsection I shall assume that an ASD null-K{\"a}hler structure 
$({\cal M}, g)$
admits
an additional parallel spinor $o_{A'}$ such that $o_{A'}\iota^{A'}=1$. 
Now there exists a covariantly
constant basis of the spin space $S^{A'}$, and $({\cal M}, g)$
is pseudo hyper-K\"ahler.
In the split signature  we can
arrange for one of the complex structures to be real and  for the other
two to be purely imaginary:
\[
-I^2=S^2=T^2=1,\qquad IST=1, 
\]
and $S$ and $T$ determine a pair of transverse null foliations.
Now \[g(X, Y)= g(IX, IY)=-g(SX, SY)=-g(TX,TY)\] for any pair of real vectors
$X, Y$. The endomorphism $I$ endows ${\cal M}$ with the structure of a 
two-dimensional complex K\"ahler manifold, as does every other complex
structure $aI+bS+cT$ parametrised by the points of the hyperboloid
$a^2-b^2-c^2=1$. Using the identification between the two-forms,
and endomorphisms given by $g$ we can write
\[
S=\Sm^{0'0'}-\Sm^{1'1'},\qquad
I=\Sm^{0'0'}+\Sm^{1'1'},\qquad T=\Sm^{0'1'}.
\]
Killing vectors  on pseudo hyper-K\"ahler spaces give rise to
a homomorphism 
\[
S^1\rightarrow SO(2,1).
\]
Therefore $K$ induces 
an $S^1$ action on the hyperboloid. Assume that $g(K, K)\neq 0$.
If the action is trivial, then
we are dealing with the $(++--)$ Gibbons-Hawking ansatz \cite{GH78}. 
Otherwise, there always 
exists a fixed point of the $S^1$ action.
If the two-form corresponding to this fixed point is non-degenerate,
then $g$ is given in terms of solutions to the $SU(\infty)$ equation
\cite{FP79, L91,W90}.
Finally if the fixed point corresponds to a degenerate  
two-form,
then $g$ is given by (\ref{HCdKP1}).
Pseudo hyper-K\"ahler metrics which admit $K$ such that $g(K, K)=0$, and
${\cal L}_Kg={\cal L}_KI={\cal L}_KS
={\cal L}_KT=0$ have been found in \cite{BG94}.
Conformal symmetries of  pseudo hyper-K\"ahler metrics with non-null
self-dual derivative have been classified in \cite{DT01}
\section{Acknowledgements}
I am grateful to David Calderbank,
Nigel Hitchin, Lionel Mason, 
Maciej Przanowski, David Robinson and Paul Tod  for useful discussions.
I would like to thank {\em Centro de Investigacion y de Estudios
Avanzados} in Mexico, where part of this work was done, for 
the hospitality.
My visit to Mexico was partially supported by grants
32427-E Proyecto de CONACYT  and LMS  5619.
\section{Appendix--spinor notation}
Let  ${\cal M}$ be a real four-manifold equipped with a
$(++--)$ metric $g$ and compatible volume form $\nu$.
We use the conventions of Penrose and Rindler \cite{PR86}.  
$a,b,...$ are
four-dimensional vector indices and $A, B, ..., A', B', ...$ are
two-dimensional spinor indices. They have ranges $0, 1$ and $0', 1'$
respectively.
The tangent space at each point of
${\cal M}$ is isomorphic to a tensor product of the two real spin spaces 
\be
\label{isomorphism}
T{\cal M}=S^A\otimes S^{A'}.
\ee
This isomorphism is given by
\[
V^a\longrightarrow V^{AA'}=
\left (
\begin{array}{cc}
V^0+V^3&V^1+V^2\\
V^1-V^2&V^0-V^3
\end{array}
\right ).
\]
Orthogonal transformations decompose into
products of ASD and SD rotations 
\be
\label{basicisom}
SO(2, 2)=(SL(2, \R)\times \widetilde{SL}(2, \R))/\Z_2.  
\ee 
The Lorentz transformation $V^a\longrightarrow {\Lambda^a}_bV^b$
is equivalent to 
\[
V^{AA'}\longrightarrow {\lambda^A}_BV^{BB'}{\lambda^{A'}}_{B'},
\]
where ${\lambda^A}_B$ and ${\lambda^{A'}}_{B'}$ are elements of $SL(2, \R)$
and $\widetilde{SL}(2, \R)$.

Spin dyads $(o^A, \iota^{A})$ and $(o^{A'}, \iota^{A'})$ span $S^A$
and $S^{A'}$ respectively.  The spin spaces $S^A$ and $S^{A'}$ are
equipped with parallel symplectic forms $\varepsilon_{AB}$ and
$\varepsilon_{A'B'}$ such that
$\varepsilon_{01}=\varepsilon_{0'1'}=1$.  These anti-symmetric objects
are used to raise and lower the spinor indices via 
$\iota_A=\iota^B\varepsilon_{BA}, \iota^A=\varepsilon^{AB}\iota_B$.
We shall use 
normalised spin frames so  that
\[
o^B\iota^C-\iota^Bo^C=\varepsilon^{BC},\;\;\;
o^{B'}\iota^{C'}-\iota^{B'}o^{C'}=\varepsilon^{B'C'}.
\] 
 Let $e^{AA'}$ be a null tetrad of one-forms on ${\cal M}$ i.e.
\[
g=\varepsilon_{AB}\varepsilon_{A'B'}e^{AA'}e^{BB'}=
2(e^{00'}e^{11'}-e^{10'}e^{01'}),
\] 
and let
$\nabla_{AA'}$ be the frame of vector fields dual to $e^{AA'}$. 
The orientation is given by fixing the volume form
\[
\nu=e^{01'}\wedge e^{10'}\wedge e^{11'}\wedge e^{00'}.
\]
Apart from orientability
${\cal M}$ must satisfy some other topological restrictions for 
$++--$ metric, and a
global spinor fields to exist. We shall not take them into
account as we work locally in ${\cal M}$.

Any two-form $\Om_{ab}$ can be written as
\[
\Om_{ABA'B'}=\Om_{AB}\varepsilon_{A'B'}+{\widetilde\Om}_{A'B'}\varepsilon_{AB}
\]
where ${\Om}_{AB}$ are  ${\widetilde\Om}_{A'B'}$ symmetric in their indices
since $\Om_{ab}$ is skew. This is the decomposition of a  two-form into its
ASD and SD part. 
The space of SD two forms is
therefore 
isomorphic to a symmetric tensor product of two primed spin spaces.

The local basis $\Sm^{AB}$ and $\Sm^{A'B'}$ of  spaces of ASD and SD
two-forms are defined by
\be\label{sdforms}
e^{AA'}\wedge e^{BB'}=\varepsilon^{AB}\Sm^{A'B'}+\varepsilon^{A'B'}\Sm^{AB}. 
\ee
The first Cartan structure equations are 
\[
\d e^{AA'}=e^{BA'}\wedge{\Gamma^{A}}_{B}+e^{AB'}\wedge{\Gamma^{A'}}_{B'},
\]
where $\Gamma_{AB}$ and $\Gamma_{A'B'}$ are the  $SL(2, \R)$
and $\widetilde{SL}(2, \R)$
spin connection one-forms.  They are symmetric in their indices, and
\[
 \Gamma_{AB}=
\Gamma_{CC'AB}e^{CC'},\;\;\Gamma_{A'B'}=\Gamma_{CC'A'B'}e^{CC'}
,\;\;\; \Gamma_{CC'A'B'}=o_{A'}\nabla_{CC'}\iota_{B'}-
\iota_{A'}\nabla_{CC'}o_{B'}.
\]
The curvature of the spin connection
\[
{R^A}_B=\d{\Gamma^A}_B+{\Gamma^A}_C\wedge{\Gamma^C}_B
\]
decomposes as
\[
{R^A}_B={C^A}_{BCD}\Sm^{CD}+(1/12)R{\Sm^{A}}_{B}+{\Phi^A}_{BC'D'}\Sm^{C'D'},
\]
and similarly  for ${R^{A'}}_{B'}$. Here $R$ is the Ricci scalar, 
$\Phi_{ABA'B'}$ is the trace-free part of the Ricci tensor $R_{ab}$,
\[
-2\Phi_{ABA'B'}=R_{ab}-\frac{1}{4}Rg_{ab},
\]
and $C_{ABCD}$ is the ASD part of the Weyl tensor
\[
C_{abcd}=\varepsilon_{A'B'}\varepsilon_{C'D'}C_{ABCD}+
\varepsilon_{AB}\varepsilon_{CD}C_{A'B'C'D'}.
\]
A conformal structure is called ASD iff $C_{A'B'C'D'}=0$.
Define the operators $\triangle_{AB}$ and $\triangle_{A'B'}$ by
\[
[\nabla_a, \nabla_b]=\varepsilon_{AB}\triangle_{A'B'}+
\varepsilon_{A'B'}\triangle_{AB}.
\]
The spinor Ricci identities are
\be
\label{ricci1}
\triangle_{AB}\iota_{A'}=\Phi_{ABA'B'}\iota^{B'}
\ee
\be
\label{ricci2}
\triangle_{A'B'}\iota_{C'}=[C_{A'B'C'D'}-\frac{1}{12}R\varepsilon_{D'(A'}
\varepsilon_{B')C'}]\iota^{D'},
\ee
(and analogous equations for unprimed spinors).
Bianchi identities translate to 
\be
\label{bianchi}
{\nabla^{A'}}_AC_{A'B'C'D'}=\nabla_{(B'}^B\Phi_{C'D')AB},\qquad
\nabla^{AA'}\Phi_{ABA'B'}+\frac{1}{8}\nabla_{BB'}R=0.
\ee

Let $K$ be a pure Killing vector. Then
$\nabla_{(A}^{(A'}K_{B)}^{B')}=0, \nabla^aK_a=0$. This implies
\[
\nabla_aK_b=\phi_{AB}\varepsilon_{A'B'}+\psi_{A'B'}\varepsilon_{AB},
\]
where $\psi_{A'B'}$ and  $\phi_{AB}$ are symmetric spinors. 
The well known identity $
\nabla_a\nabla_bK_c=R_{bcad}K^d $ implies
\be
\label{Killing_identity}
{\nabla^A}_{A'}\psi_{B'C'}=-2C^{D'}_{A'B'C'}K^{A}_{D'}-
2{K^B}_{(A'}\Phi_{B'C')B}^A+\frac{R}{6}\varepsilon_{A'(B'}{K^A}_{C')}
-\frac{4}{3}\varepsilon_{A'(B'}\Phi_{C')}^{D'DA}K_{DD'}.
\ee
Therefore in an ASD vacuum $\psi_{A'B'}=const$.
A Lie derivative of a spinor along a Killing vector $K$ is given by
\[
{\cal L}_K{\iota^{A'}}=K^b\nabla_b\iota^{A'}-\psi^{A'}_{B'}\iota^{B'}.
\]
Note that if $\iota_{A'}$ is covariantly constant, and   
the Killing vector $K$ preserves the $\iota_{A'}$ 
we deduce that  $\psi_{A'B'}=\psi\iota_{A'}\iota_{B'}$ for some function
$\psi$. The identity (\ref{Killing_identity}) yields
\[
\nabla_{AA'}\psi=\frac{1}{2}F_{AB}K^{B}_{A'}
\]
so $\psi\neq const$ unless $g$ is hyper-K\"ahler. 

\end{document}